\documentclass[11pt]{article}
\usepackage{amssymb,amsfonts,amsmath,latexsym,epsf,tikz,url}

\newtheorem{theorem}{Theorem}[section]

\newtheorem{definition}[theorem]{Definition}

\newcommand{\proof}{\noindent{\bf Proof. }}
\newcommand{\qed}{\hfill $\square$\medskip}

\begin{document}

\title{The distinguishing number and the distinguishing  index of line and graphoidal graphs}

\author{
Saeid Alikhani
\and
Samaneh Soltani  $^{}$\footnote{Corresponding author}
}

\date{\today}

\maketitle

\begin{center}
Department of Mathematics, Yazd University, 89195-741, Yazd, Iran\\
{\tt alikhani@yazd.ac.ir, s.soltani1979@gmail.com}
\end{center}


\begin{abstract}
The distinguishing number (index) $D(G)$ ($D'(G)$) of a graph $G$ is the least integer $d$
such that $G$ has an vertex labeling (edge labeling)  with $d$ labels  that is preserved only by a trivial automorphism. 
 A graphoidal cover of $G$ is a collection $\psi$ of (not necessarily open) paths in $G$ such that every path in $\psi$ has at least two vertices, every vertex of $G$ is an internal vertex of at most one path in $\psi$ and every edge of $G$ is in exactly one path in $\psi$. Let $\Omega(G,\psi)$  denote the intersection graph of $\psi$. 
 A graph $H$ is called a graphoidal graph, if there exists a graph $G$ and
a graphoidal cover $\psi$ of $G$ such that $H\cong \Omega (G, \psi)$. In this paper, we study the distinguishing number and the distinguishing  index of the line graph and the graphoidal graph of a simple connected graph $G$.
\end{abstract}


\medskip
\noindent{\bf AMS Subj.\ Class.}: 05C25 

\section{Introduction and definitions}
Let $G = (V,E)$ be a simple, finite, connected and undirected graph, and let ${\rm Aut}(G)$ be its  \textit{automorphism group}.
     A labeling of $G$, $\phi : V \rightarrow \{1, 2, \ldots , r\}$, is said to be \textit{$r$-distinguishing},  if no non-trivial  automorphism of $G$ preserves all of the vertex labels. The point of the labels on the vertices is to destroy the symmetries of the
graph, that is, to make the automorphism group of the labeled graph trivial.
Formally, $\phi$ is $r$-distinguishing if for every non-identity $\sigma \in {\rm Aut}(G)$, there exists $x$ in $V$ such that $\phi(x) \neq \phi(\sigma(x))$. The \textit{distinguishing number} of a graph $G$ is defined  by
\begin{equation*}
D(G) = {\rm min}\{r \vert ~ G ~\text{{\rm has a labeling that is $r$-distinguishing}}\}.
\end{equation*} 

This number has defined in \cite{Albert}. Similar to this definition, in \cite{Kali1} the \textit{distinguishing index} $D'(G)$ of $G$, which is  the least integer $d$
such that $G$ has an edge colouring   with $d$ colours that is preserved only by a trivial
automorphism, has been  defined. $D'(G)$ can be arbitrary smaller than $D(G)$, for example $D'(K_{p,p})=2$ and $D(K_{p,p})=p+1$, for $p\geq 4$.

The \textit{line graph}  $L(G)$ of a graph $G$ is  the graph whose vertices are edges of $G$ and  two edges $e, e' \in   V (L(G)) = E(G)$ are adjacent if they share an endpoint in common.  
The line graphs can be  viewed a special case of   graphoidal graphs.  The concept of graphoidal cover was introduced by Acharya and Sampathkumar \cite{Acharya}.  
\begin{definition}{\rm \cite{Acharya}}
 A graphoidal cover of a graph $G$ is a collection $\psi$ of (not necessarily open) paths in $G$ satisfying the following conditions:
\begin{itemize}
\item[(i)] Every path in $\psi$ has at least two vertices.
\item[(ii)] Every vertex of $G$ is an internal vertex of at most one path in $\psi$.
\item[(iii)]  Every edge of $G$ is in exactly one path in $\psi$.
\end{itemize}
\end{definition}
\begin{definition}{\rm \cite{Acharya}}
If $F = \{S_1, S_2,\ldots , S_n\}$ is  a family of distinct nonempty
subsets of a set $S$, whose union is $S$, then the intersection graph of $F$, denoted by
$\Omega­(F)$, is the graph whose vertex set and edge   set are given by
\begin{equation*}
V­(F) =  \{S_1, S_2,\ldots , S_n\},~~E­(F) = \{S_iS_j :~ i \neq j ,~ S_i \cap S_j \neq \emptyset \}.
\end{equation*}
\end{definition}

Let $\mathcal{G}_G$ be the set of all graphoidal covers of $G$ and $\psi \in  \mathcal{G}_G$. The intersection graph on $\psi$ is denoted by $\Omega­(G, \psi)$. 
 A graph $H$ is called a \textit{graphoidal graph} if there exists a graph $G$
and $\psi  \in  \mathcal{G}_G$ such that $H \cong \Omega ­(G, \psi)$.
Since $E(G)$ is obviously a graphoidal cover of $G$, all line graphs are graphoidal graphs. In \cite{Acharya} Acharya and Sampathkumar have proved that all the Beineke's forbidden
subgraphs  of line graphs (see Figure \ref{exceptionssubgraphs}) are graphoidal graphs and hence they conjectured that all graphs
are graphoidal. In \cite{Arumugam} Arumugam and Pakkiam disproved this conjecture by showing that complete
bipartite graphs $K_{m, n}$ with $mn> 2(m + n)$ are not graphoidal graphs. They also obtained a
forbidden subgraph characterization of all bipartite graphs which are graphoidal.
Note that if $P= (v_0, v_1, \ldots , v_k)$ is a path, not necessarily
open, in a graph $G$, then $v_0$ and $v_k$ are called \textit{terminal vertices} 
and $v_1, \ldots , v_{k-1}$ are called \textit{internal vertices}  of $P$. For cycles
(considered as closed paths), there is an inherent “ordering”
of vertices as in paths. So, when we say that a cycle $C$ of a
graph $G$ is a member of a graphoidal cover $\psi$ of $G$, we should mention the
vertex at which the cycle $C$ begins, and this particular vertex
is considered as the terminal vertex of $C$ and all other vertices
on $C$ are called internal vertices of $C$. Given a graphoidal cover $\psi$ of $G$, a
vertex $v$  is said to be interior to $\psi$ if $v$ is an internal vertex of
an element of $\psi$ and is called exterior to $\psi$ otherwise.

\medskip
In the next section, we study the distinguishing number and the distinguishing index of line graphs. In Section 3, we consider graphoidal graphs and study their distinguishing number and distinguishing index. 

\section{Study of $D(G)$ and $D'(G)$ for the line graphs}
\begin{figure}
	\begin{center}
		\includegraphics[width=0.67\textwidth]{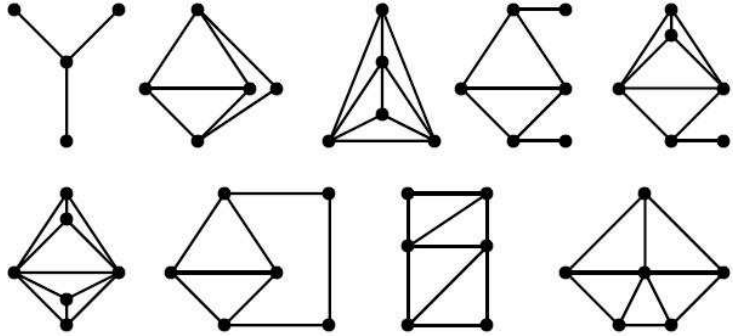}
		\caption{\label{exceptionssubgraphs} The nine  forbidden subgraph characterization of line graphs.}
	\end{center}
\end{figure}

For a simple graph $G$,  the line graph  $L(G)$ is  the graph whose vertices are edges of $G$ and  two edges $e, e' \in   V (L(G)) = E(G)$ are adjacent if they share an endpoint in common.  The  following theorem characterize the line graphs.
\begin{theorem}{\rm \cite{Beineke}}\label{characterizinglinegraphs}
The following statements are equivalent for a graph $G$.
\begin{enumerate}
\item[(i)] The graph $G$ is the line graph of some graphs.
\item[(ii)]  The edges of $G$ can be partitioned into complete subgraphs in such a
way that no vertex belongs to more than two of the subgraphs.
\item[(iii)]  The graph $K_{1,3}$ is not an induced subgraph of $G$; and if $abc$ and $bcd$ are
distinct odd triangles, then $a$ and $d$ are adjacent.
\item[(iv)]  None of the nine graphs in Figure \ref{exceptionssubgraphs} is an induced subgraph of $G$.
\end{enumerate}
\end{theorem}

\medskip
To study the distinguishing number and the distinguishing  index of  $L(G)$, we need more information about the automorphism group of $L(G)$.
 Let  $\gamma_G : {\rm Aut} (G) \rightarrow {\rm Aut} (L(G))$ be given by $(\gamma_G \phi)(\{u, v\}) = \{\phi(u), \phi(v)\}$ for every $\{u, v\} \in  E(G)$. In \cite{Sabidussi}, Sabidussi proved the following Theorem which we need it later.
\begin{theorem}{\rm \cite{Sabidussi}}\label{autlinegraph}
 Suppose that $G$ is a connected graph that is not $P_2, Q$, or $L(Q)$ (see Figure \ref{autline}). Then 
$G$ is a group isomorphism, and so ${\rm Aut}(G) \cong {\rm Aut}(L(G))$.
\begin{figure}
	\begin{center}
		\includegraphics[width=0.55\textwidth]{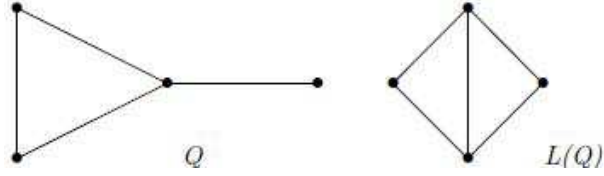}
		\caption{\label{autline} graphs $Q$ and $L(Q)$ of Theorem \ref{autlinegraph}.}
	\end{center}
\end{figure}
\end{theorem}

Now we are ready to obtain the distinguishing number of line graph of a graph. We note that if $G= L(Q)$, then it is easy to see that $D'(L(Q))=2$, while $D(L(L(Q)))=3$. 
\begin{theorem}\label{disnumline}
Suppose that $G$ is a connected graph that is not $P_2$ and $L(Q)$. Then $D(L(G))= D'(G)$.
\end{theorem}
\proof  If $G = Q$, then it is easy to see that $D'(Q) = D(L(Q))=2$.  If $G \neq Q$,  First we show that $D(L(G))\leq D'(G)$. Let $c:E(G)\rightarrow \{1, \ldots , D'(G)\}$ be an edge distinguishing labeling of $G$. We define $c': V(L(G)) \rightarrow \{1, \ldots , D'(G)\}$ such that $c'(e) = c(e)$, where $e\in V(L(G))  = E(G)$. The vertex labeling $c'$ is a distinguishing vertex labeling of $L(G)$, because if $f$ is an automorphism of $L(G)$ preserving the labeling, then $c'(f(e))= c'(e)$, and hence $c(f(e))= c(e)$ for any $e\in E(G)$. On the other hand,  by Theorem \ref{autlinegraph}, $f = \gamma_G \phi$ for some automorphism $\phi$ of $G$. Thus from  $c(f(e))= c(e)$ for any $e\in E(G)$, we can conclude that $c(\gamma_G \phi (e))= c(e)$ and so $c(\{\phi(u), \phi(v)\})= c(\{u,v\})$ for every $\{u,v\}\in E(G)$. This means that $\phi$ is an automorphism of $G$ preserving the labeling $c$, and so $\phi$ is the identity automorphism of $G$. Therefore $f$ is the identity automorphism of $L(G)$, and hence $D(L(G))\leq D'(G)$.  For the converse, suppose that  $c:V(L(G))\rightarrow \{1, \ldots , D(L(G))\}$ is a vertex distinguishing labeling of $L(G$). We define $c': E(G) \rightarrow \{1, \ldots ,D(L(G))\}$ such that $c'(e) = c(e)$ where $e\in    E(G)$. The edge labeling $c'$ is a distinguishing edge labeling of $G$. Because if $f$ is an automorphism of $G$ preserving the labeling, then $c'(f(e))= c'(e)$, and hence $c(f(e))= c(e)$ for any $e\in E(G)$. Then, there exist the automorphism $\gamma_G f$ of $L(G)$ such that  $\gamma_G f (\{u,v\})= \{f(u), f(v)\}$ for every $\{u,v\}\in E(G)$, by Theorem \ref{autlinegraph}. Thus from  $c(f(e))= c(e)$ for any $e\in E(G)$, we can conclude that  $c(\{u,v\})= c(\{f(u), f(v)\})= c (\gamma_G f (\{u,v\}))$ for every $\{u,v\}\in E(G)$, which means that  $\gamma_G f$ preserves the distinguishing  vertex labeling of $L(G)$, and hence $\gamma_G f$ is the identity automorphism of $L(G)$. Therefore $f$ is the identity automorphism of $G$, and so $D'(G) \leq D(L(G))$. \qed

An upper bound for the distinguishing index of line graphs   follows immediately from the following result of Pil\'sniak in \cite{nord}. A $K_{1,3}$-free graph, called also a \textit{claw-free graph}, is a graph containing no copy of $K_{1,3}$ as an induced subgraph. 
\begin{theorem}{\rm \cite{nord}}\label{clawfree}
If $G$ is a connected claw-free graph, then $D'(G)\leq  3$.
\end{theorem}
Now by Part (iii) of Theorem \ref{characterizinglinegraphs} and Theorem \ref{clawfree}, we have the following theorem.
\begin{theorem} 
If $G$ is a connected graph, then $D'(L(G))\leq  3$.
\end{theorem}

\section{Study of $D(G)$ and $D'(G)$ for graphoidal graphs}

\begin{figure}
	\begin{center}
		\includegraphics[width=0.8\textwidth]{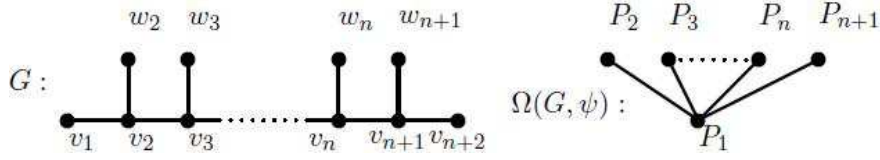}
		\caption{\label{dd} Graph $G$ in the proof of Theorem \ref{lar}.}
	\end{center}
\end{figure}

We recall that  a graphoidal cover of a graph $G$ is a collection $\psi$ of (not necessarily open) paths in $G$ satisfying the following conditions: every path in $\psi$ has at least two vertices,   every vertex of $G$ is an internal vertex of at most one path in $\psi$, and  every edge of $G$ is in exactly one path in $\psi$.
Let $\psi$ be a graphoidal cover of $G$ and $\Omega(G,\psi)$ denote the intersection graph of $\psi$. Thus the vertices of $\Omega(G,\psi)$ are the paths in $\psi$ and two paths in $\psi$ are adjacent in $\Omega(G,\psi)$  if and only if they have a common vertex. A graph $G$ is said to be graphoidal if there exists a graph $H$ and a graphoidal cover $\psi$ of $H$ such that $G$ is isomorphic to   $\Omega(H,\psi)$. 
First we state and prove the following theorem: 
\begin{theorem} \label{lar} 
	There exists a graph $G$ such that the values $|D(G)-D(\Omega (G, \psi))|$ and $|D'(G)-D'(\Omega (G,\psi))|$    can be arbitrarily large.
	\end{theorem} 
	\proof
	 Let $G$ be a graph as shown in Figure \ref{dd}. It is easy to see that $D(G) = D'(G)=2$. A graphoidal cover of $G$ is $\psi=\{P_1, P_2,\ldots , P_{n+1}\}$ where $P_1 = (v_1,v_2, \ldots , v_{n+2})$ and $P_i = (v_i, w_i)$ for $2\leq i \leq n+1$. Thus $\Omega (G, \psi)= K_{1,n}$, 
	 and hence $D(\Omega (G, \psi))=D'(\Omega (G, \psi))=n$. Therefore we have the result.\qed

\medskip
To study the distinguishing index of graphoidal graphs, we need some theorems. 
\begin{theorem}{\rm{\cite{Kali1}}}\label{boundindex}
	If $G$ is a connected graph of order $n \geq 3$, then $D'(G) \leq \Delta(G)$, 
	unless $G$ is $C_3, C_4$ or $C_5$.
\end{theorem}
\begin{theorem}{\rm{\cite{nord}}}\label{boundindex2}
	Let $G$ be a connected graph that is neither a symmetric nor a bisymmetric
	tree. If the maximum degree of $G$ is at least 3, then
	$D'(G) \leq \Delta(G) - 1$, unless $G$ is $K_4$ or $K_{3,3}$.
\end{theorem}
\begin{theorem}{\rm{\cite{Collins}}}\label{treebound1}
	If $T$ is a tree of order $n \geq 3$, then $D(T ) \leq \Delta(T)$. Furthermore, equality is achieved if
	and only if $T$ is a symmetric tree or a path of odd length.
\end{theorem}
Before we state the next theorem we need to define some terms:
Let $v$ be a vertex of a graph $G$ and let $f$ be a $k$-distinguishing edge labeling of $G$.
We say that f is \textit{$v$-distinguishing} if $f$ is  preserved only by the trivial automorphism
among the automorphisms $\alpha \in Aut(G)$ for which $\alpha(v) = v$ holds. If $vw$ is an
edge of a tree $T$, then let $T_v$ and $T_w$ be the components of $T - \{vw\}$ with $v \in T_v$ and $w \in T_w$. A family $\mathcal{T}$ consists of those trees $T$ of order at least 3, for which the
following conditions are fulfilled: (1) $T$ is a bicentric tree with the central edge $e = vw$, (2) $T_v$ and $T_w$ are isomorphic trees, (3) $T_v$ admits a unique $v$-distinguishing edge $D(T)$-labeling.
\begin{theorem}{\rm{\cite{alikhani&sandi}}}\label{treebound2}
	Let $T$ be a tree of order $n \geq 3$. Then $D'(T) = D(T ) + 1$ if $T$ belongs to $\mathcal{T}$, and $D'(T ) = D(T )$ for all other trees.
\end{theorem}

\begin{figure}
	\begin{center}
		\includegraphics[width=0.9\textwidth]{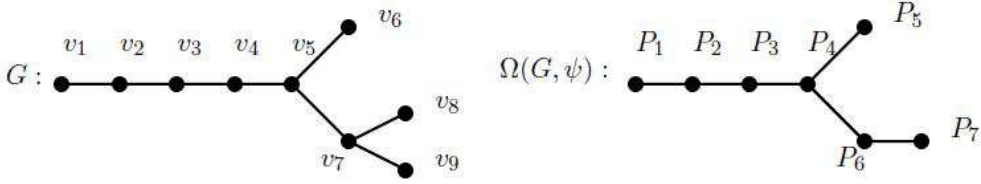}
		\caption{\label{cc} An example for sharpness of inequality of Corollary \ref{distnumbgraphoidalpath}.}
	\end{center}
\end{figure}

 Now, we obtain bounds for  the distinguishing index of graphoidal graphs. 
\begin{theorem}
Let $G$ be a connected graph and $\psi$ be a graphoidal cover of $G$ such that the order of $\Omega(G, \psi)$ is at least $3$. If  $\Omega(G, \psi)\neq C_3, C_5$, then $1\leq D'(\Omega(G, \psi))\leq |\psi|-1$. Moreover, $D'(\Omega(G, \psi))= |\psi|-1$ if and only if $\Omega(G, \psi)$ is $C_4$, $K_4$ or $K_{1, |\psi|-1}$.
\end{theorem}
\proof By Figure \ref{cc}, we can conclude the left inequality. For the right inequality,  by contradiction, we suppose that $D'(\Omega(G, \psi))\geq |\psi|$. Since $D'(\Omega(G, \psi))\leq \Delta (\Omega(G, \psi))$ by Theorem \ref{boundindex}, we get $|\psi|\leq \Delta (\Omega(G, \psi))$, which is a contradiction, since $\Delta (\Omega(G, \psi)) \leq |V(\Omega(G, \psi))|-1= |\psi|-1$. 

 If $D'(\Omega(G, \psi))= |\psi|-1$, then $D'(\Omega(G, \psi))= \Delta(\Omega(G, \psi))$, by Theorem \ref{boundindex}, unless for $\Omega(G, \psi) = C_4$. If $\Delta(\Omega(G, \psi))\geq 3$, then  we have  $D'(\Omega(G, \psi))= \Delta(\Omega(G, \psi))$, if and only if  $\Omega(G, \psi)$ is  a symmetric, a bisymmetric tree,  $K_4$ or $K_{3,3}$,  by Theorem \ref{boundindex2}. Now since $\Delta(\Omega(G, \psi))=|\psi|-1$, so $D'(\Omega(G, \psi))= |\psi|-1$ if and only if $\Omega(G, \psi)$ is $C_4$, $K_4$ or $K_{1, |\psi|-1}$, by Theorems \ref{treebound1} and \ref{treebound2}.\qed

\medskip

\begin{figure}
	\begin{center}
		\includegraphics[width=0.94\textwidth]{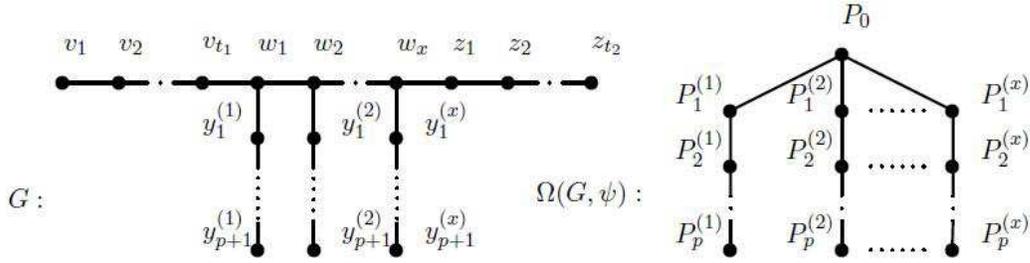}
		\caption{\label{ee} The graph $G$  and its graphoidal of Theorem \ref{exsitence}. }
	\end{center}
\end{figure}

The following theorem states that for every $i$, $2\leq i \leq |\psi|-1$, there exists a connected graphoidal graph with the distinguishing index $i$. 

\begin{theorem}\label{exsitence}
Suppose that  $x$ and $p$ are positive integers and $(x,p)\neq (1,1)$. If $k = 1+xp$ and $d = \lceil \sqrt[p]{x}\rceil$, then there exist a graph $G$ of order $n \geq x(p+2)$ and graphoidal cover $\psi$ of $G$ with $|\psi|=k$ such that $D'(\Omega (G, \psi))=d$.
\end{theorem}
\proof  Let $G$ be a graph in Figure \ref{ee}.  If $\psi =\{P_0, P_1^{(1)}, \ldots , P_p^{(1)}, \ldots , P_1^{(x)}, \ldots , P_p^{(x)}\}$ where $P_0 = (v_1, \ldots , v_{t_1}, w_1, \ldots , w_x, z_1, \ldots , z_{t_2})$ and for every $1 \leq i \leq x$,  $P_1^{(i)} = (w_i,y_1^{(i)})$,   $P_j^{(i)} = (y_j^{(i)},y_{j+1}^{(i)})$ for each $2\leq j \leq p$, then $\Omega(G, \psi)$ is as shown in Figure \ref{ee}. It is easy to compute that $D'(\Omega(G, \psi))= \lceil \sqrt[p]{x} \rceil$. With respect to integers  $t_1$ and $t_2$, we get that $n \geq x(p+2)$.\qed

\medskip
In sequel, we gives bounds for the distinguishing number of the graphoidal graphs. 
\begin{theorem}\label{distnumbgraphoidal}
Let $G$ be a connected graph of order $n \geq 3$. Then $$D'(G)-2\leq D(\Omega (G, \psi))\leq |\psi|.$$

\end{theorem}
\proof Since $\Omega (G, \psi)$ has $|\psi|$ vertices, so it is clear that $D(\Omega (G, \psi))\leq |\psi|$. To prove the left inequality,  let $D(\Omega (G, \psi)) = t$ and $X_i = \{P_{i1}, \ldots , P_{is_i}\}$, $1 \leq i \leq t$ be the set of vertices of $\Omega (G, \psi)$ having label $i$ in the distinguishing labeling of $\Omega (G, \psi)$.  Now we want to label the edges of $G$ distinguishingly, using   $t+2$ labels. For this purpose,   for every $i$, $1 \leq i \leq t$,  we label the edge set of  all paths of  $\psi$  (not necessarily open) in $X_i$ of length $j$ with $j$-tuple  $(i, t+1, \underbrace{t+2,\ldots , t+2}_{(j-2)-times})$ of labels. We note that since there exist the closed paths (cycles)  in $\psi$, so we need at most three different labels to distinguishing  these closed paths.  Then we have an edge labeling of $G$ with $t+2$ labels. To show that this edge labeling is distinguishing, we prove that if $f$  is an automorphism of $G$ preserving this edge labeling, then $f$ maps $\psi$ to $\psi$, setwise. In fact, if $P\in \psi$, then $f(P)$ is a path of the same length as $P$. Since $f$ preserves  the edge labeling of $G$,  so the label of edges  of $f(P)$ is the same as edges of $P$. With respect to the method of lebeling of edge set of each  path in $\psi$, we conclude that  $f(P) \in \psi$.  Hence $f$ maps  a path in $\psi$ to a path in $\psi$.  Thus $f$ can be considered  as an automorphism of  $\Omega (G, \psi)$ preserves the distinguishing vertex labeling of  $\Omega (G, \psi)$ with $t$ labels. Thus $f$ is the identity automorphism of  $\Omega (G, \psi)$, and this means that $f(P) = P$ for all $P\in \psi$. On the other hand, since we labeled  the edge set of each path in $\psi$ distinguishingly with at most three labels, so $f$ fixes  all vertices of path $P$, where  $P\in \psi$. Thus $f$ is the identity automorphism of $G$.   \qed

The  bounds of Theorem \ref{distnumbgraphoidal} are sharp. For the left inequality, it is sufficient to consider $G=C_i$, $3\leq i \leq 5$, and $\psi = \{C_i\}$. Thus  $\Omega (G, \psi)= K_1$, and hence $D'(G)=3$ and $ D(\Omega (G, \psi))=1$. To show that the upper bound of  this theorem is sharp, we consider $G=C_3$ and  $\psi = \{e_1,e_2, e_3\}$. Thus  $\Omega (G, \psi)= C_3$, and hence $|\psi |=3$ and $ D(\Omega (G, \psi))=3$.

\begin{theorem}\label{distnumbgraphoidalpath}
Let $G$ be a connected graph of order $n \geq 3$. If graphoidal cover $\psi$ of $G$ contains only open paths, then $D'(G)\leq D(\Omega (G, \psi))+1$.
\end{theorem}
\proof The proof is exactly the same as proof of Theorem \ref{distnumbgraphoidal}, except the method of labeling  of edges of paths in $\psi$. Since all paths in $\psi$ are open, then by notation of Theorem \ref{distnumbgraphoidal}, for every $i$, $1 \leq i \leq t$,   we label the edge set of all open paths of  $\psi$   in $X_i$ of length $j$ with the  $j$-tuple    $(i, \underbrace{t+1,\ldots , t+1}_{(j-1)-times})$ of labels. Then we have an edge labeling of $G$ with $t+1$ labels. To show  this edge labeling is distinguishing, we prove that if $f$  is an automorphism of $G$ preserving this edge labeling, then $f$ maps $\psi$ to $\psi$ setwise. In fact, if $P\in \psi$, then $f(P)$ is a path of the same length as $P$. Since $f$ preserves  the edge labeling of $G$,  so the label of edges  of $f(P)$ is the same as edges of $P$. With respect to the method of lebeling of edge set of each  path in $\psi$, we conclude that  $f(P) \in \psi$.  Hence $f$ maps  a path in $\psi$ to a path in $\psi$.  Thus $f$ can be considered  as an automorphism of  $\Omega (G, \psi)$ preserves the distinguishing vertex labeling of  $\Omega (G, \psi)$ with $t$ labels. Thus $f$ is the identity automorphism of  $\Omega (G, \psi)$, and this means that $f(P) = P$ for all $P\in \psi$. On the other hand, since we labeled  the edge set of each path in $\psi$ distinguishingly with at most two labels, so $f$ fixes  all vertices of path $P$ where  $P\in \psi$. Thus $f$ is the identity automorphism of $G$.  \qed

The bound of Theorem \ref{distnumbgraphoidalpath} is sharp.   let $G$ be as shown in Figure \ref{cc}. It is easy  to obtain that $D'(G)=2$. A graphoidal cover of $G$ is $\psi=\{P_1, P_2,\ldots , P_7\}$ where $P_1 = (v_1,v_2)$, $P_2 = (v_2,v_3)$, $P_3 = (v_3,v_4)$, $P_4 = (v_4,v_5)$, $P_5 = (v_5,v_6)$, $P_6 = (v_5,v_7, v_8)$, and $P_7 = (v_8,v_9)$. Thus $\Omega (G, \psi)$ is an asymmetric graph, and hence $D(\Omega (G, \psi))=1$. 

\medskip
We conclude the paper with the following theorem, which shows that the value $|D'(G)-D(\Omega(G, \psi))|$ can be arbitrary.
\begin{theorem}
For every $i$, $ i\geq 0$, there exist a connected graph $G$ and  a graphoidal cover $\psi$ of $G$ such that  $|D'(G)-D(\Omega(G, \psi))|=i$.
\end{theorem}
\proof  Let $G$ and $\Omega(G, \psi)$ be the two graphs of Theorem \ref{exsitence}. Since the automorphism group of $G$ has at most  one nonidentity automorphism, so $D'(G)\leq 2$. On the other hand $\Omega (G, \psi)$ is a tree with a central vertex, and so $D(\Omega (G, \psi))= D'(\Omega (G, \psi))= \lceil \sqrt[p]{x}\rceil$, by Theorem \ref{treebound2}. Now, with respect to the values of $x$ and $p$, we obtain the result.\qed

\end{document}